\documentclass[12pt,oneside]{amsart}
\usepackage{amssymb,latexsym,amsmath,amsfonts}
\usepackage{graphicx}

\newcommand{\Ff}{{\mathbb F}}

\newcommand{\p}{{\mathbf P}}

\newcommand{\A}{{\mathbf A}}

\newcommand{\s}{{\mathbf S}}

\newcommand{\Proof}{ \noindent{\bf Proof:}\quad }

\def\QED{\qed\medskip\par}
\newtheorem{Theorem} {Theorem} [section]

\newtheorem{Lemma} [Theorem] {Lemma}

\def\Sq{\;\Box}
\def\Nsq{\;\not\!\!\Box}

\begin{document}

\title{Some $p$-ranks related to a conic in $PG(2,q)$}

\author{Junhua Wu}
\address{Department of Mathematical Sciences, Worcester Polytechnic Institute, Worcester, MA, 01609 USA} \email{wuj@wpi.edu}


\begin{abstract}
Let $\A$ be the incidence matrix of lines and points of the classical projective plane $PG(2,q)$ with $q$ odd. With respect to a conic in $PG(2,q)$, the matrix $\A$ is partitioned into $9$ submatrices. The rank of each of these submatices over $\Ff_q$, the defining field of $PG(2,q)$, is determined.
\end{abstract}

\maketitle

\section{Introduction}
Let $\Ff_q$ be the finite field of order $q$, where $q=p^e$, $p$ is
an odd prime, and $e$ is a positive integer. Let $PG(2,q)$ denote the
classical projective plane. That is, the points and lines of
$PG(2,q)$ are the $1$-dimensional subspaces and $2$-dimensional
subspaces of $\Ff_q^3$, respectively, and the incidence is the
natural inclusion. An {\it oval} in $PG(2,q)$ is a set of $q+1$
points, no three of which are collinear. A {\it conic} in $PG(2,q)$
is the set of points $\langle (x,y,z)\rangle$ satisfying a nonzero
quadratic form. A conic is said to be {\it nondegenrate} if it does
not contain an entire line of $PG(2,q)$. It is well known
\cite{hirsch} that every nondegenerate conic is an oval; moreover,
by a linear change of coordinates, any nondegenerate conic is
equivalent to
$$\mathcal{O} = \{\langle(1,t,t^2)\rangle\mid t\in \Ff_q\}\cup\{\langle(0,0,1)\rangle\},$$ the
set of projective solutions of the nondegenerate quadratic form
$$Q(X_0,X_1,X_2) = X_1^2-X_0 X_2$$ over $\Ff_q$.

In the case where $p$ is odd, Segre \cite{seg} proved that an oval
in $PG(2,q)$ must be a nondegenerate conic. Therefore in $PG(2,q)$,
$q$ an odd prime power, ovals and nondegenerate conics are
essentially the same objects. In the rest of this paper, we will
always assume that $p$ is an odd prime, and use the above
$\mathcal{O}$ as our ``standard" conic. With respect to
$\mathcal{O}$, a line $\ell$ is called {\it secant}, {\it tangent},
or {\it skew} according as $|\ell\cap\mathcal{O}|=2$, $1$, or $0$.
We use $Se$, $T$, and $Sk$ to denote the set of secant lines, the
set of tangent lines, and the set of skew lines, respectively. Since
$\mathcal{O}$ is an oval, every line of $PG(2,q)$ must fall into one
of these sets. A point of $PG(2,q)$ is called {\it external}, {\it
absolute}, or {\it internal} with respect to $\mathcal{O}$ according
as it lies on 2, 1, or 0 tangent lines to $\mathcal{O}$. We use $E$
and $I$ to denote the set of external points and the set internal
points, respectively.

Let $\A$ be the line-point incidence matrix of $PG(2,q)$. That is,
the rows of $\A$ are labeled by the lines of $PG(2,q)$, the columns of
$\A$ are labeled by the points of $PG(2,q)$, and the $(\ell, P)$-
entry $\A$ is 1 if $P\in \ell$, $0$ otherwise. The following result
is a special case of the well-known $p$-rank formula for the
incidence matrix of the Singer design (see \cite{smith},
\cite{MacMann}, and \cite{GDel}).
\begin{Theorem}\label{rankofA}
rank$_p(\A)$ = $\binom{p+1}{2}^e+1$.

\end{Theorem}

In a recent paper \cite{DMM}, Droms, Mellinger and Meyer considered
the following partition of $\A$ into submatrices:
\begin{displaymath}
\A=\left(\begin{array}{ccc}
\A_{11} & \A_{12} & \A_{13} \\
\A_{21} & \A_{22} & \A_{23} \\
\A_{31} & \A_{32} & \A_{33} \\
\end{array}\right),
\end{displaymath}
where the rows of $\A_{11}, \A_{21}$ and $\A_{31}$ are labeled by
the tangent, skew, and secant lines, respectively, and the columns
of $\A_{11}, \A_{12}$ and $\A_{13}$ are labeled by the absolute,
internal, and external points, respectively. These authors used the
submatrices $A_{ij}$ to construct binary linear codes. Some of these
codes are good examples of structured low-density parity-check
(LDPC) codes, cf \cite{DMM}. Based on computational evidence, the
authors of \cite{DMM} made conjectures on the dimensions of the
binary LDPC codes. These conjectures were investigated in
\cite{thesis}, and will be proved in forthcoming papers. In this
paper, we are interested in the $p$-ranks of the submatrices
$\A_{ij}$, where $p$ is the characteristic of $\Ff_q$.

We mention that several authors have considered the $p$-ranks of
various submatrices of $\A$. Let $\A_{nonsec}$ (resp. $\A_{sec}$)
denote the submatrices of $\A$ obtained by deleting the rows of $\A$
indexed by the secant (resp. tangent and skew) lines. That is,
\begin{displaymath}
\A_{nonsec} =\left(\begin{array}{ccc}
\A_{11} & \A_{12} & \A_{13} \\
\A_{21} & \A_{22} & \A_{23} \\
\end{array}\right)\;\text{and}\;
\A_{sec}=\left(\begin{array}{ccc}
\A_{31} & \A_{32} & \A_{33}\\
\end{array}\right).
\end{displaymath}
Then Blokhuis and Moorhouse \cite{blokhuis} proved the following
theorem:
\begin{Theorem}\label{sub_1}
Use the above notation, we have
\begin{itemize}
\item[(i)] $rank_p \A_{nonsec}$ = $\binom{p+1}{2}^e+1$,
\item[(ii)] $rank_p \A_{sec}$ = $\binom{p+1}{2}^e$.
\end{itemize}
\end{Theorem}
The authors of \cite{blokhuis} also computed the 2-ranks of
$\A_{nonsec}$ and $\A_{sec}$ when $q$ is a power of 2. Carpenter
\cite{carp}, and later Kamiya and Fossorier \cite{kf} computed the
2-ranks of some more refined submatrices of $\A$ when $q$ is a power
of 2. Our main result in this paper is the following theorem:
\begin{Theorem}\label{you}
Let $\A_{ij}$, $1\le i, j \le 3$, be the submatrices of $\A$ defined
above. Then
\begin{itemize}
\item[(i)] rank$_p(\A_{11})$ = $q+1$,
\item[(ii)] rank$_p(\A_{12})$ = rank$_p(\A_{21})$ = $0$,
\item[(iii)] rank$_p(\A_{13})$ = rank$_p(\A_{31})$ = $q$,
\item[(iv)] rank$_p(\A_{22}) = $rank$_p(\A_{23})$ = rank$_p(\A_{32})$ = $\binom{p+1}{2}^e - q$,
\item[(v)] rank$_p(\A_{33})$ = $\binom{p+1}{2}^e$.
\end{itemize}
\end{Theorem}

The paper is organized as follows. In Section 2, we collect basic
geometric results related to a conic in $PG(2,q)$. In Section 3, we
first construct several vector spaces of polynomials over $\Ff_q$;
we then convert the problem of computing the $p$-ranks of the above
submatrices into the one of calculating the dimensions of these polynomial
spaces over $\Ff_q$. Finally we prove Theorem~\ref{you} by
explicitly computing the dimensions of the polynomial spaces
constructed. The proofs mainly rely on the following
Nullstellensatz proved in \cite{blokhuis}:

\begin{Lemma}\label{new}
Let $F_d[X_0,X_1,X_2]$ be the vector space of homogeneous
polynomials of degree $d$ in $\Ff_q[X_0,X_1,X_2]$, together with 0.
Let $f(X_0,X_1,X_2)\in F_d[X_0,X_1,X_2]$, with $d\leq q-1$ and $q$
odd. Define $Q(X_0,X_1,X_2)=X_1^2-X_0X_2$.
\begin{itemize}
\item[(i)] If $f(x_0,x_1,x_2) = 0$ whenever $Q(x_0,x_1,x_2)$ is a non-zero square, where $(x_0,x_1,x_2)\in \Ff_q^3$, then $f = 0$.
\item[(ii)] If $f(x_0,x_1,x_2) = 0$ whenever $Q(x_0,x_1,x_2)$ is a non-square or zero, where $(x_0,x_1,x_2)\in \Ff_q^3$, then $f=0$.
\end{itemize}
\end{Lemma}

\section{Some known geometric results related to a conic in $PG(2,q)$}

A {\it correlation} of $PG(2,q)$ is a bijection sending points to
lines and lines to points that reverses inclusion. A {\it polarity}
of $PG(2,q)$ is a correlation of order $2$. The image of a point
$\p$ under a correlation $\sigma$ is denoted by  $\p^\sigma$, and
that of a line $\ell$ is denoted by $\ell^\sigma$. It can be shown
\cite[p.~181]{hirsch} that the nondegenerate quadratic form
$Q(X_0,X_1,X_2)$ = $X_1^2-X_0X_2$ induces a polarity $\sigma$ of
$PG(2,q)$, which can be represented by the matrix
\begin{displaymath}
M=\left(\begin{array}{ccc}
0 & 0 & -\frac{1}{2} \\
0 & 1 & 0 \\
-\frac{1}{2} & 0 & 0 \\
\end{array}\right).
\end{displaymath}
For example, if $\p=\langle(x,y,z)\rangle$ is a point of $PG(2,q)$,
then its image under $\sigma$ is $\p^{\sigma}=\langle
M(x,y,z)^{\top}\rangle=\langle(z,-2y,x)^{\top}\rangle$.

\begin{Lemma}\label{bijection}$(${\rm{\cite[p.~181]{hirsch}}}$)$
The polarity $\sigma$ defines the following three bijections:
$$\sigma:\;I\rightarrow\;Sk,$$
$$\sigma:\;E\rightarrow\;Se,$$
and
$$\sigma:\;\mathcal{O}\rightarrow\;T.$$
\end{Lemma}
For convenience, we denote the set of all non-zero squares of
$\Ff_q$ by $\Sq_q$, and the set of non-squares by $\Nsq_q$. Then we
have the following lemma.

\begin{Lemma}\label{lem1}$(${\rm{\cite[p.~182]{hirsch}}}$)$
A line $\langle(r,m,n)^{\top}\rangle$ of $PG(2,q)$ is skew, tangent,
or secant to $\mathcal{O}$ if and only if $m^2-4nr \in \Nsq_q$,
$m^2-4nr = 0$, or $m^2-4nr\in \Sq_q$, respectively.
\end{Lemma}

Using the polarity $\sigma$ induced by $Q$ and the previous lemma,
we have:

\begin{Lemma}\label{lem5}$(${\rm{\cite[p.~182]{hirsch}}}$)$
A point $\p=\langle(x,y,z)\rangle$ of $PG(2,q)$ is internal,
external, or absolute if and only if $y^2-xz \in \Nsq_q$, $y^2-xz
\in \Sq_q$, or $y^2-xz = 0$, respectively.
\end{Lemma}
\begin{Lemma}{\rm (\cite[p.~170]{hirsch})}
We have the following tables:
\begin{table}[htp]
\begin{center}
\caption{Number of lines of various types}
\begin{tabular}{ccc}
\hline
{Tangent lines} & {Skew lines} & {Secant lines} \\
\hline
{$q+1$} & $\frac{q(q-1)}{2}$ & $\frac{q(q+1)}{2}$ \\
\hline
\end{tabular}
\label{table:points}
\end{center}
\end{table}

\begin{table}[htp]
\begin{center}
\caption{Number of points of various types}
\begin{tabular}{ccc}
\hline
{Absolute points} & {Internal points} & {External points}\\
\hline
$q+1$& $\frac{q(q-1)}{2}$ & $\frac{q(q+1)}{2}$ \\
\hline
\end{tabular}
\label{table:lines}
\end{center}
\end{table}

\begin{table}[htp]
\begin{center}
\caption{Number of points on lines of various types}
\bigskip
\begin{tabular}{cccc}
\hline
{Name} & {Absolute points} & {External points} & {Internal points} \\
\hline
{Tangent lines} & $1$ & $q$ & $0$ \\
{Secant lines} & $2$ & $\frac{q-1}{2}$ & $\frac{q-1}{2}$ \\
{Skew lines} & $0$ & $\frac{q+1}{2}$ & $\frac{q+1}{2}$\\
\hline
\end{tabular}
\label{table:tab1}
\end{center}
\end{table}

\begin{table}[htp]
\begin{center}
\caption{Number of lines through points of various types}
\bigskip
\begin{tabular}{cccc}
\hline
{Name} & {Tangent lines} & {Secant lines} & {Skew lines} \\
\hline
{Absolute points} & $1$ & $q$ & $0$ \\
{External points} & $2$ & $\frac{q-1}{2}$ & $\frac{q-1}{2}$ \\
{Internal points} & $0$ & $\frac{q+1}{2}$ & $\frac{q+1}{2}$\\
\hline
\end{tabular}
\label{tab2}
\end{center}
\end{table}
\end{Lemma}

\section{Some Vector Spaces of Polynomials}
Let $\p=\langle(x_0,x_1,x_2)\rangle$ be a point of $PG(2,q)$ and let
$\ell=\langle(y_0,y_1,y_2)^{\top}\rangle$ be a line of $PG(2,q)$. We
have $\p\in \ell$ if and only if $x_0y_0+x_1y_1+x_2y_2 = 0$. So the
$(\ell,\p)$-entry of $\A$ is given by
\begin{displaymath}
(\A)_{\ell,\p} = 1 - (x_0y_0+x_1y_1+x_2y_2)^{q-1} =\begin{cases}
1, & \text{if} \; x_0y_0+x_1y_1+x_2y_2 =0,\\
0, & \text{otherwise}.

\end{cases}
\end{displaymath}

Now we define a $(q^3-1)\times q^3$ matrix $\mathbf{S}$, whose rows
are indexed by vectors ${\bf
y}=(y_0,y_1,y_2)\in\Ff_q^3\setminus\{(0,0,0)\}$, whose columns are
indexed by vectors ${\bf x}=(x_0,x_1,x_2)\in \Ff_q^3$, and the
$({\bf y},{\bf x})$-entry is equal to $1 -
(x_0y_0+x_1y_1+x_2y_2)^{q-1}$. Similar to the partition of $\A$
defined in Section 1, we partition $\mathbf{S}$ into the following
form
\begin{displaymath}
\mathbf{S}=\left(
\begin{array}{ccc}
\mathbf{S}_{11}& \mathbf{S}_{12} & \mathbf{S}_{13} \\
\mathbf{S}_{21}& \mathbf{S}_{22} & \mathbf{S}_{23} \\
\mathbf{S}_{31}& \mathbf{S}_{32} & \mathbf{S}_{33} \\
\end{array}
\right),
\end{displaymath}
where the rows of $\mathbf{S}_{11}$, $\mathbf{S}_{21}$, and
$\mathbf{S}_{31}$ are labeled by vectors $(y_0,y_1,y_2)\in
\Ff_q^3\setminus\{(0,0,0)\}$ such that $y_1^2-4y_0y_2=0$,
$y_1^2-4y_0y_2\in\Nsq_q$, and $y_1^2-4y_0y_2\in\Sq_q$, respectively,
and the columns of $\mathbf{S}_{11}$, $\mathbf{S}_{12}$, and
$\mathbf{S}_{13}$ are labeled by vectors $(x_0,x_1,x_2)\in \Ff_q^3$
such that $x_1^2-x_0x_2=0$, $x_1^2-x_0x_2\in\Nsq_q$, and
$x_1^2-x_0x_2\in\Sq_q$, respectively. We also define four more
submatrices of $\mathbf{S}$:
\begin{displaymath}
\mathbf{S}_{nonsec} =\left(\begin{array}{ccc}
\mathbf{S}_{11} & \mathbf{S}_{12} & \mathbf{S}_{13} \\
\mathbf{S}_{21} & \mathbf{S}_{22} & \mathbf{S}_{23} \\
\end{array}\right),\;
\mathbf{S}_{sec}=\left(\begin{array}{ccc} \mathbf{S}_{31} &
\mathbf{S}_{32} & \mathbf{S}_{33}\end{array}\right),
\end{displaymath}
\begin{displaymath}
\mathbf{S}_{sk}=\left(\begin{array}{ccc} \mathbf{S}_{21} &
\mathbf{S}_{22} & \mathbf{S}_{23}\end{array}\right),\;\text{and}\;
\mathbf{S}_{T}=\left(\begin{array}{ccc}\mathbf{S}_{11}& 
\mathbf{S}_{12} & \mathbf{S}_{13} \end{array}\right).
\end{displaymath}

\begin{Lemma}\label{equal}
Use the above notation, we have
\begin{itemize}
\item[(i)] $rank_p\A = rank_p\mathbf{S}$.
\item[(ii)] $rank_p\mathbf{S}_{i j} = rank_p\mathbf{A}_{i j}$ for $1\le i, j \le 3$, except $(i,j)=(2,1)$.
\item[(iii)] $rank_p\A_{nonsec} = rank_p \mathbf{S}_{nonsec}$, $rank_p \A_{sec} = rank_p \mathbf{S}_{sec}$.
\item[(iv)] $rank_p\mathbf{S}_{sk}=rank_p\A_{sk}$.
\end{itemize}
\end{Lemma}

{\Proof} Note that the rows of $\s$ are indexed by the vectors
$(y_0,y_1,y_2)\in\Ff_q^3\setminus\{(0,0,0)\}$ whose 
transposes represent the lines of $PG(2,q)$. Assume that the first
column of $S$ is indexed by $(0,0,0)$. Hence each entry of the first 
column of $\s$ is $1$. Since each line contains 
$q+1$ points and each point can be represented by
 $q-1$ different non-zero vectors, 
we see that the sum of all the columns of $\s$ is a $q^2$ 
multiple of an all one column vector of the proper size, 
which is a zero column vector over $\Ff_q$. This implies that 
the first column of $\s$ is a linearly combination of all the other 
columns of $\s$ over $\Ff_q$. Hence, the matrix $\s^{'}$ obtained 
by deleting the first column of $\s$ has the same $p$-rank as $\s$. 
By deleting duplicate rows and columns of $\s^{'}$ and 
permuting the rows and columns of the resulting matrix, 
we can obtain the matrix $\A$. This shows that $\s$, $\s^{'}$, and $\A$ 
have the same $p$-rank. So $(i)$ follows.

It is clear that $rank_p\s_{ij} = rank_p\A_{ij}$ for $(i,j)=(1,3)$, $(2,2)$, $(2,3)$, $(3,2)$, 
and $(3,3)$. Note that each tangent line contains a point of $\mathcal{O}$ 
and each point can be represented by $q-1$ non-zero vectors. 
So the sum of all the columns of $S_{11}$ is a zero column over $\Ff_q$, 
which indicates that the first column of $\s_{11}$ 
is a linear combination of all other columns of $\s_{11}$. 
Thus, by deleting the first column of $\s_{11}$ and 
the duplications of rows and columns of $\s_{11}$ and 
permuting the rows and columns of the resulting matrix, 
we get $\A_{11}$. Hence, $rank_p\s_{11} = rank_p\A_{11}$. 
Similarly, the first column of $\s_{31}$ is a linear combination of 
all the other columns of $\s_{31}$ since, again, the sum of all the columns 
of $\s_{31}$ is a zero column by noting that each secant line contains 
$2$ points of $\mathcal{O}$. Thus, by deleting the first column of $\s_{31}$ 
and the duplications of rows and columns of $\s_{31}$ and 
permuting the rows and columns of the resulting matrix, 
we get $\A_{31}$. Hence, $rank_p\s_{31} = rank_p\A_{31}$. 
It is also clear that $rank_p\s_{21}=1$ and $rank_p\A_{21} =0$. So $(ii)$ is proved.

The proofs of $(iii)$ and $(iv)$ are essentially the same as the proofs of 
$(i)$ and $(ii)$. We omit the detail.
\QED

\begin{Lemma}\label{in_1}
$rank_p \mathbf{S}_{sk} = \binom{p+1}{2}^e-q$.
\end{Lemma}
{\Proof} Note that the rows of $\A_{nonsec}$ indexed by the tangent lines are linear independent since $\A_{11}$ is a permutation matrix. When deleting the rows of $\A_{nonsec}$ indexed by the tangent lines, the $p$-rank of the matrix decreases by $q+1$. This implies that $rank_p\A_{sk}=\binom{p+1}{2}^e-q$ by $(i)$ of Theorem~\ref{sub_1}. From $(iv)$ of Lemma~\ref{equal}, the lemma follows.
\QED

Now we define
\begin{displaymath}
\begin{array}{lll}
\mathcal{Z}(Sk) & =&\{(y_0,y_1,y_2)\in \Ff_q^3\mid y_1^2-4y_0y_2\in \Nsq_q\},\\
\mathcal{Z}(Se) & =&\{(y_0,y_1,y_2)\in \Ff_q^3\mid y_1^2-4y_0y_2\in \Sq_q\},\\
\mathcal{Z}(T) &=&\{(y_0,y_1,y_2)\in \Ff_q^3\mid y_1^2-4y_0y_2 = 0\}\setminus\{(0,0,0)\}.\\
\end{array}
\end{displaymath}
Also
\begin{displaymath}
\begin{array}{lll}
\mathcal{Z}(I) &=&\{(x_0,x_1,x_2)\in \Ff_q^3\mid
x_1^2-x_0x_2\in\Nsq_q\},\\
\mathcal{Z}(E) &=&\{(x_0,x_1,x_2)\in \Ff_q^3\mid x_1^2-x_0x_2\in
\Sq_q\},\\
\mathcal{Z}(\mathcal{O})&=&\{(x_0,x_1,x_2)\in \Ff_q^3\mid x_1^2-
x_0x_2
= 0\}.\\
\end{array}
\end{displaymath}

Let $\Ff_q[X_0,X_1,X_2]$ be the vector space of polynomials in three
indeterminates over $\Ff_q$ and $F_{q-1}[X_0,X_1,X_2]$ be the
subspace of $\Ff_q[X_0,X_1,X_2]$ over $\Ff_q$ consisting of the
homogeneous polynomials of degree $q-1$, together with $0$. Define
the following subspaces of $\Ff_q[X_0,X_1,X_2]$ over $\Ff_q$:
\begin{displaymath}
\begin{array}{llll}
\mathcal{M}&=& \left\langle 1 - (y_0X_0+y_1X_1 + y_2X_2)^{q-1}\mid  (y_0,y_1,y_2)\in \Ff_q^3\right\rangle,\\
\mathcal{M}^{Sk} &=& \left\langle1 - (y_0X_0+y_1X_1 + y_2X_2)^{q-1}\mid (y_0,y_1,y_2)\in \mathcal{Z}(Sk)\right\rangle,\\

\mathcal{M}^{Se}&=& \left\langle1 - (y_0X_0+y_1X_1 + y_2X_2)^{q-1}\mid (y_0,y_1,y_2)\in \mathcal{Z}(Se)\right\rangle,\\

\mathcal{M}^{T}&= &\left\langle1 - (y_0X_0+y_1X_1 + y_2X_2)^{q-1}\mid (y_0,y_1,y_2)\in \mathcal{Z}(T)\right\rangle.\\
\end{array}
\end{displaymath}

Given a vector space $M$ over $\Ff_q$, we use $dim_q M$ to denote
its dimension over $\Ff_q$. If $\mathbf{B}$ is a matrix over $\Ff_q$, then
we use row$(\mathbf{B})$ to denote the span of the rows of $\mathbf{B}$
over $\Ff_q$.

\begin{Lemma}\label{dim_1} Use the above notation,
\begin{itemize}
\item[(i)] $dim_q \mathcal{M}^{Sk} = \binom{p+1}{2}^e-q$,
\item[(ii)] $dim_q\mathcal{M}^{Se} = \binom{p+1}{2}^e$,
\item[(iii)] $dim_q \mathcal{M}^{T} = q+1$.
\end{itemize}
\end{Lemma}
{\Proof}  Define
$$\psi_0: \mathcal{M}^{Sk}\rightarrow \text{row}({\mathbf{S}_{sk}})$$
by
$$g(X_0,X_1,X_2)\mapsto \left(g(x_0,x_1,x_2): (x_0,x_1,x_2)\in \Ff_q^3\right).$$
Note that $\psi_0([1-(y_0X_0+y_1X_1+y_2X_2)^{q-1}])$ is simply the row
of $\mathbf{S}_{sk}$ indexed by $(y_0,y_1,y_2)$ with
$y_1^2-4y_0y_2\in \Nsq_q$. Let $$g(X_0,X_1,X_2)
=\displaystyle\sum_{(y_0,y_1,y_2)\in
\mathcal{Z}(Sk)}a_{y_0,y_1,y_2}\left[1-(y_0X_0+y_1X_1+y_2X_2)^{q-1}\right]
\in \mathcal{M}^{Sk}.$$ Then
\begin{displaymath}
\begin{array}{lll}
\psi_0(g(X_0,X_1,X_2)) & = & \displaystyle\sum_{(y_0,y_1,y_2)\in \mathcal{Z}(sk)}a_{y_0,y_1,y_2}\mathbf{S}_{y_0,y_1,y_2}^{sk}\\
{}& = & \displaystyle\sum_{(y_0,y_1,y_2)\in \mathcal{Z}(sk)}a_{y_0,y_1,y_2}\psi_0([1-(y_0X_0+y_1X_1+y_2X_2)^{q-1}])\in\text{row}(\mathbf{S}_{sk}),\\
\end{array}
\end{displaymath}
where $\mathbf{S}_{y_0,y_1,y_2}^{sk}$ is the row vector of
$\mathbf{S}_{sk}$ indexed by $(y_0,y_1,y_2)$. The first equality
above comes from the definition of $\psi_0$. So $\psi_0$ is a surjective
$\Ff_q$-linear map. Now assume that
$g(X_0,X_1,X_2)\in\text{Ker}(\psi_0)$. Then $g(x_0,x_1,x_2) = 0 $ for
all $(x_0,x_1,x_2)\in \Ff_q^3$. So $g(X_0,X_1,X_2) = 0$. 
Thus, we proved that
$$\mathcal{M}^{Sk}\cong \text{row}(\mathbf{S}_{sk})$$ as vector
spaces over $\Ff_q$. Part $(i)$ of the lemma follows immediately from
Lemma~\ref{in_1}.

The proofs of part $(ii)$ and $(iii)$ are similar. We need only check that the maps
$$\psi_1: \mathcal{M}^{Se}\rightarrow \text{row}({\mathbf{S}_{se}})$$
defined by
$$g(X_0,X_1,X_2)\mapsto \left(g(x_0,x_1,x_2): (x_0,x_1,x_2)\in \Ff_q^3\right)$$
and
$$\psi_2: \mathcal{M}^{T}\rightarrow \text{row}({\mathbf{S}_{T}})$$
defined by
$$g(X_0,X_1,X_2)\mapsto \left(g(x_0,x_1,x_2): (x_0,x_1,x_2)\in \Ff_q^3\right)$$
are both $\Ff_q$-isomorphisms. Since they are the same as the proof of $(i)$, we omit the details.
\QED

\begin{Theorem}\label{in_2}
As $\Ff_q$-spaces,
\begin{itemize}
\item[(i)]$\mathcal{M}^{Sk}\cong \text{row}(\mathbf{S}_{22})$,
\item[(ii)]$\mathcal{M}^{Se}\cong \text{row}(\mathbf{S}_{33})$,
\item[(iii)] $\mathcal{M}^{Sk}\cong \text{row}(\mathbf{S}_{23})$,
\item[(iv)] $\mathcal{M}^{T}/\langle\mathbf{J}(X_0,X_1,X_2)\rangle\cong \text{row}(\mathbf{S}_{13})$, where $\langle\mathbf{J}(X_0,X_1,X_2)\rangle$ is the $1$-dimensional subspace generated by the polynomial $$\mathbf{J}(X_0,X_1,X_2) = 1+\displaystyle\sum_{(y_0,y_1,y_2)\in\mathcal{Z}(T)}(y_0X_0+y_1X_1+x_2Y_2)^{q-1}$$ over $\Ff_q$.
\end{itemize}
\end{Theorem}
{\Proof} Define the map
$$\gamma: \mathcal{M}^{sk}\rightarrow \text{row}(\mathbf{S}_{22})$$
by
$$g(X_0,X_1,X_2)\mapsto \left(g(x_0,x_1,x_2): (x_0,x_1,x_2)\in \mathcal{Z}(I)\right).$$ We can show that $\gamma$ is a well-defined surjective $\Ff_q$-linear map by applying arguments similar to the ones of Lemma~\ref{dim_1}. Now let $h(X_0,X_1,X_2)\in\text{Ker}(\gamma)$. Then

\begin{equation}
\begin{array}{lll}
h(x_0,x_1,x_2) &= & \displaystyle\sum_{(y_0,y_1,y_2)\in \mathcal{Z}(Sk)}a_{y_0,y_1,y_2}\left[1-(x_0y_0+x_1y_1+x_2y_2)^{q-1}\right]\\
{} & = & \displaystyle\sum_{(y_0,y_1,y_2)\in\mathcal{Z}(Sk)}a_{y_0,y_1,y_2} - \displaystyle\sum_{(y_0,y_1,y_2)\in\mathcal{Z}(Sk)}a_{y_0,y_1,y_2}(x_0y_0+x_1y_1+x_2y_2)^{q-1}\\
{}& = & 0
\end{array}
\end{equation}
for all $(x_0,x_1,x_2)\in\mathcal{Z}(I)$, where $a_{y_0,y_1,y_2}\in \Ff_q$. Note that $(0,0,0)\notin\mathcal{Z}(I)$.

Let $\p=\langle(x_0,x_1,x_2)\rangle$ be any internal point and Sk$_{\p}$ be the set of $\frac{q+1}{2}$ skew lines through $\p$.  We denote the set of row vectors whose transposes represent the skew lines through $\p$ by $\mathcal{Z}(Sk_{\p})$. Then $(1)$ gives
\begin{displaymath}
h(x_0,x_1,x_2)  =  \displaystyle\sum_{(y_0,y_1,y_2)\in\mathcal{Z}(Sk_{\p})} a_{y_0,y_1,y_2} = 0.
\end{displaymath}
Hence
\begin{equation}
\displaystyle\sum_{\p=\langle(x_0,x_1,x_2)\rangle\in I}h(x_0,x_1,x_2) = \displaystyle\sum_{\p=\langle(x_0,x_1,x_2)\rangle\in I}\displaystyle\sum_{(y_0,y_1,y_2)\in\mathcal{Z}(Sk_{\p})}a_{y_0,y_1,y_2} = 0.
\end{equation}
Since both $\langle(y_0,y_1,y_2)^{\top}\rangle$ and $\langle(x_0,x_1,x_2)\rangle$ are represented by $(q-1)$ different non-zero vectors, we see that $|\mathcal{Z}(Sk_{\p})|=\frac{(q+1)(q-1)}{2}$ by the last row of Table~\ref{tab2} and $|\mathcal{Z}(I)|=\frac{q(q-1)^2}{2}$ by the Table~\ref{table:points}. So $(2)$ can be written as
\begin{displaymath}
\begin{array}{lll}
\displaystyle\sum_{\p=\langle(x_0,x_1,x_2)\rangle\in I}\displaystyle\sum_{(y_0,y_1,y_2)\in\mathcal{Z}(Sk_{\p})}a_{y_0,y_1,y_2} &=& \frac{q+1}{2}\cdot\left(\displaystyle\sum_{(y_0,y_1,y_2)\in\mathcal{Z}(Sk)}a_{y_0,y_1,y_2}\right) \\
{} &=& 0.
\end{array}
\end{displaymath}
As $p\nmid\frac{q+1}{2}$ for odd $p$, we must have
\begin{displaymath}
\begin{array}{lll}
\displaystyle\sum_{(y_0,y_1,y_2)\in\mathcal{Z}(Sk)}a_{y_0,y_1,y_2}& = & 0.
\end{array}
\end{displaymath}
Thus $h(X_0,X_1,X_2)\in F_{q-1}[X_0,X_1,X_2]$ and
$$h(x_0,x_1,x_2) = \displaystyle\sum_{(y_0,y_1,y_2)\in\mathcal{Z}(Sk)}a_{y_0,y_1,y_2}(x_0y_0+x_1y_1+x_2y_2)^{q-1} = 0 $$
for each $(x_0,x_1,x_2)\in\mathcal{Z}(I)$. Hence $h(X_0,X_1,X_2) = 0$ by $(ii)$ of Lemma~\ref{new}. We have proved that $\gamma$ is an $\Ff_q$-isomorphism. So $(i)$ follows.

Part $(ii)$ can be proved in the same fashion. Consider the map
$$\beta: \mathcal{M}^{Se}\rightarrow \text{row}(\mathbf{S}_{33})$$
by
$$g(X_0,X_1,X_2)\mapsto \left(g(x_0,x_1,x_2): (x_0,x_1,x_2)\in \mathcal{Z}(E)\right).$$
It is easy to see that $\beta$ is a well-defined surjective $\Ff_q$-linear map.  Let $h(X_0,X_1,X_2)\in\text{Ker}(\beta)$. Then

\begin{equation}
\begin{array}{lll}
h(x_0,x_1,x_2) &= & \displaystyle\sum_{(y_0,y_1,y_2)\in \mathcal{Z}(Se)}a_{y_0,y_1,y_2}\left[1-(x_0y_0+x_1y_1+x_2y_2)^{q-1}\right]\\
{} & = & \displaystyle\sum_{(y_0,y_1,y_2)\in\mathcal{Z}(Se)}a_{y_0,y_1,y_2} - \displaystyle\sum_{(y_0,y_1,y_2)\in\mathcal{Z}(Se)}a_{y_0,y_1,y_2}(x_0y_0+x_1y_1+x_2y_2)^{q-1}\\
{}& = & 0
\end{array}
\end{equation}
for all $(x_0,x_1,x_2)\in\mathcal{Z}(E)$, where $a_{y_0,y_1,y_2}\in \Ff_q$. Note that $(0,0,0)\notin\mathcal{Z}(E)$.

Let $\p=\langle(x_0,x_1,x_2)\rangle$ be any external point and Se$_{\p}$ be the set of secant lines through $\p$. We denote the set of different vectors $(y_0,y_1,y_2)$ whose transposes represent the secant lines through $\p$ by $\mathcal{Z}(Se_{\p})$. Then $(3)$ gives
\begin{displaymath}
h(x_0,x_1,x_2)  =  \displaystyle\sum_{(y_0,y_1,y_2)\in\mathcal{Z}(Se_{\p})} a_{y_0,y_1,y_2} = 0.
\end{displaymath}
Hence
\begin{equation}\label{sum_2}
\displaystyle\sum_{\p=\langle(x_0,x_1,x_2)\rangle\in E}h(x_0,x_1,x_2) = \displaystyle\sum_{\p=\langle(x_0,x_1,x_2)\rangle\in E}\displaystyle\sum_{(y_0,y_1,y_2)\in\mathcal{Z}(Se_{\p})}a_{y_0,y_1,y_2} = 0.
\end{equation}
Since $|\mathcal{Z}(Se_{\p})|=\frac{(q-1)^2}{2}$ by the last row of Table~\ref{tab2} and $|\mathcal{Z}(E)|=\frac{q(q+1)(q-1)}{2}$ by the Table~\ref{table:points}, we see that $(\ref{sum_2})$ can be written as
\begin{displaymath}
\begin{array}{llll}
\displaystyle\sum_{\p=\langle(x_0,x_1,x_2)\rangle\in E}\displaystyle\sum_{(y_0,y_1,y_2)\in\mathcal{Z}(Se_{\p})}a_{y_0,y_1,y_2} &= &\frac{q-1}{2}\cdot\left(\displaystyle\sum_{(y_0,y_1,y_2)\in\mathcal{Z}(Se)}a_{y_0,y_1,y_2}\right)\\
{} & = & 0.\\
\end{array}
\end{displaymath}
As $p\nmid\frac{q-1}{2}$ for odd $p$, we must have
\begin{displaymath}
\begin{array}{lll}
\displaystyle\sum_{(y_0,y_1,y_2)\in\mathcal{Z}(Se)}a_{y_0,y_1,y_2} & = & 0.
\end{array}
\end{displaymath}
Thus $h(X_0,X_1,X_2) \in F_{q-1}[X_0,X_1,X_2]$ and
\begin{displaymath}
\begin{array}{lll}
h(x_0,x_1,x_2)& = &\displaystyle\sum_{(y_0,y_1,y_2)\in\mathcal{Z}(Se)}a_{y_0,y_1,y_2}(x_0y_0+x_1y_1+x_2y_2)^{q-1} \\
{} & = & 0.
\end{array}
\end{displaymath}
for all $(x_0,x_1,x_2)\in\mathcal{Z}(E)$. Hence $h(X_0,X_1,X_2) = 0$ by $(i)$ of Lemma~\ref{new}. So $(ii)$ is proved.

The proof of $(iii)$ is essentially the same as the proof of $(i)$, so we omit the details.

To prove $(iv)$, consider the map
$$\eta: \mathcal{M}^{T}\rightarrow \text{row}(\mathbf{S}_{13})$$
by
$$g(X_0,X_1,X_2)\mapsto \left(g(x_0,x_1,x_2): (x_0,x_1,x_2)\in \mathcal{Z}(E)\right).$$
Again, it is easy to see that $\eta$ is a well-defined surjective $\Ff_q$-linear map. Let $h(X_0,X_1,X_2)\in\text{Ker}(\eta)$. Then
\begin{equation}
\begin{array}{lll}
h(x_0,x_1,x_2) &= & \displaystyle\sum_{(y_0,y_1,y_2)\in \mathcal{Z}(T)}a_{y_0,y_1,y_2}\left[1-(x_0y_0+x_1y_1+x_2y_2)^{q-1}\right]\\
{} & = & \displaystyle\sum_{(y_0,y_1,y_2)\in\mathcal{Z}(T)}a_{y_0,y_1,y_2} - \displaystyle\sum_{(y_0,y_1,y_2)\in\mathcal{Z}(T)}a_{y_0,y_1,y_2}(x_0y_0+x_1y_1+x_2y_2)^{q-1}\\
{}& = & 0
\end{array}
\end{equation}
for all $(x_0,x_1,x_2)\in\mathcal{Z}(E)$, where $a_{y_0,y_1,y_2}\in \Ff_q$. 
Note that $(0,0,0)\notin\mathcal{Z}(E)$.

Let $\p=\langle(x_0,x_1,x_2)\rangle$ be any external point. We use $\mathcal{Z}(T_{\p})$ to denote the row vectors whose transposes represent either of two tangent lines through $\p$. Then $(5)$ simplifies to
\begin{displaymath}
\begin{array}{lll}
h(x_0,x_1,x_2) &=& \displaystyle\sum_{(y_0,y_1,y_2)\in \mathcal{Z}(T_{\p})}a_{y_0,y_1,y_2}\\
{}& =& 0.\\
\end{array}
\end{displaymath}
Let $E_{\ell}$ denote the set of external points on a given tangent line $\ell$ and $\langle\ell\rangle$ be the $1$-dimensional subspace over $\Ff_q$ generated by the row vector whose transpose represents $\ell$. Since $T_{\p}$ consists of the two tangent lines through $\p$ (the second row of Table~\ref{tab2}), then, in the multiset $\bigcup_{\p\in E_{\ell}}T_{\p}$, each tangent line other than $\ell$ appears exactly once, and $\ell$ appears exactly $q$ times. Hence,
\begin{displaymath}
\begin{array}{llll}
\displaystyle\sum_{\p=\langle(x_0,x_1,x_2)\rangle\in E_{\ell}}h(x_0,x_1,x_2) &= &\displaystyle\sum_{\p=\langle(x_0,x_1,x_2)\rangle\in E_{\ell}}\displaystyle\sum_{(y_0,y_1,y_2)\in \mathcal{Z}(T_{\p})}a_{y_0,y_1,y_2} \\
{} & = & q\cdot\left(\displaystyle\sum_{(y_0,y_1,y_2)\in\langle\ell\rangle}a_{y_0,y_1,y_2}\right) + \displaystyle\sum_{(y_0,y_1,y_2)\in\mathcal{Z}(T)\setminus\langle\ell\rangle}a_{y_0,y_1,y_2}\\
{} & = & \displaystyle\sum_{(y_0,y_1,y_2)\in\mathcal{Z}(T)\setminus\langle\ell\rangle}a_{y_0,y_1,y_2} \\
{} & = & 0.
\end{array}
\end{displaymath}
Hence for any two different tangent lines $\ell_1$ and $\ell_2$, we have
\begin{displaymath}
\begin{array}{lll}
\displaystyle\sum_{(y_0,y_1,y_2)\in\mathcal{Z}(T)\setminus\langle\ell_1\rangle}a_{y_0,y_1,y_2} &=&  \displaystyle\sum_{(y_0,y_1,y_2)\in\mathcal{Z}(T)\setminus\langle\ell_2\rangle}a_{y_0,y_1,y_2},
\end{array}
\end{displaymath}
which indicates that
\begin{displaymath}
\begin{array}{lll}
\displaystyle\sum_{(y_0,y_1,y_2)\in \langle\ell_1\rangle} a_{y_0,y_1,y_2} &= & \displaystyle\sum_{(y_0,y_1,y_2)\in \langle\ell_2\rangle} a_{y_0,y_1,y_2} = L\in \Ff_q\\
\end{array}
\end{displaymath}
for any two tangent lines $\ell_1$ and $\ell_2$. Thus,
\begin{equation}\label{eq6}
\begin{array}{llllll}
h(X_0,X_1,X_2) & = &  \displaystyle\sum_{(y_0,y_1,y_2)\in\mathcal{Z}(T)}a_{y_0,y_1,y_2} - \displaystyle\sum_{(y_0,y_1,y_2)\in\mathcal{Z}(T)}a_{y_0,y_1,y_2}(X_0y_0+X_1y_1+X_2y_2)^{q-1}\\
{}& = & \displaystyle\sum_{\ell\in T}\displaystyle\sum_{(y_0,y_1,y_2)\in\langle\ell\rangle}a_{y_0,y_1,y_2} - \displaystyle\sum_{\ell\in T}\left(\displaystyle\sum_{(y_0,y_1,y_2)\in \langle\ell\rangle}a_{y_0,y_1,y_2}\right)(y_0X_0+y_1X_1+y_2X_2)^{q-1}\\
{} & = & L\cdot(q+1) - L\cdot\left[\displaystyle\sum_{\langle(y_0,y_1,y_2)^T\rangle\in T}(X_0y_0+X_1y_1+X_2y_2)^{q-1}\right]\\
{}&=& L + L\cdot\left[(q-1)\cdot\displaystyle\sum_{\langle(y_0,y_1,y_2)^T\rangle\in T}(X_0y_0+X_1y_1+X_2y_2)^{q-1}\right]\\
{} & = & L\cdot\left[ 1 + \displaystyle\sum_{(y_0,y_1,y_2)\in \mathcal{Z}(T)}(X_0y_0+X_1y_1+X_2y_2)^{q-1}\right].\\

\end{array}
\end{equation}
The second and the last equalites in (\ref{eq6}) hold since for any two vectors $(y_0,y_1,y_2),\;(y_0^{'},y_1^{'},y_2^{'})\in\langle\ell\rangle$, $(X_0y_0+X_1y_1+X_2y_2)^{q-1}=(X_0y_0^{'}+X_1y_1^{'}+X_2y_2^{'})^{q-1}$. So $(iv)$ follows immediately.
\QED

\begin{Lemma}\label{you1}
Let $\mathbf{S}_{ij}$'s with $1\le i, j \le 3$ be the submatrices of $\mathbf{S}$. Then
\begin{itemize}
\item[(i)] $rank_p(\mathbf{S}_{11}) = q+1$,
\item[(ii)] $rank_p(\mathbf{S}_{12}) = 0$,
\item[(iii)] $rank_p(\mathbf{S}_{13}) = rank_p(\mathbf{S}_{31}) = q$,
\item[(iv)] $rank_p(\mathbf{S}_{22}) =  rank_p(\mathbf{S}_{23}) = rank_p(\mathbf{S}_{32}) = \binom{p+1}{2}^e - q$,
\item[(v)] $rank_p(\mathbf{S}_{33}) = \binom{p+1}{2}^e$.
\end{itemize}

\end{Lemma}
{\Proof} $(i)$ and $(ii)$ are clear. By $(iv)$ of Theorem~\ref{in_2}, we have $rank_p(\mathbf{S}_{13})=q$. Since $\mathbf{S}_{31}$ can be obtained by permuting the rows and columns of $\mathbf{S}_{13}$, then $rank_p(\mathbf{S}_{31}) = q$. $(iv)$ follows from $(i)$ and $(iii)$ of Theorem~\ref{in_2} and $(ii)$ of Lemma~\ref{dim_1}. $(iv)$ follows from $(ii)$ of Theorem~\ref{in_2} and the fact that $\s_{32}$ can be obtained by permuting the rows and columns of $\s_{23}$. Finally, $(v)$ follows from $(ii)$ of Lemma~\ref{dim_1} and $(ii)$ of Theorem~\ref{in_2}.
\QED

\noindent{{\bf Proof of Theorem~\ref{you}:}} The theorem follows immediately from Corollary~\ref{equal} and Lemma~\ref{you1}. \QED

\end{document}